\documentclass[11pt]{amsart}
\usepackage{amsmath}
\usepackage{mathtools}
\usepackage{wasysym}
\usepackage{MnSymbol}
\usepackage{thmtools}											
\usepackage[letterpaper,margin=1.25in]{geometry}    			 
\usepackage[english]{babel}										
\usepackage[pdfencoding=auto, psdextra, draft=false]{hyperref}
\usepackage{bookmark}
\usepackage{url}												
\usepackage[T1]{fontenc}
\usepackage{xspace}		
\usepackage{fancyhdr}
\usepackage{enumerate}
\usepackage{mathrsfs}
\usepackage{graphicx}
\usepackage{soul,color}
\usepackage{mathtools}
\usepackage{tikz-cd}
\usepackage[style=alphabetic]{biblatex}
\usepackage{csquotes}
\usepackage{chngcntr}
\counterwithin{equation}{section}

\newcommand{\horrule}[1]{\rule{\linewidth}{#1}}

\newcommand{\Title}[1]{
  \title[#1]{
    \horrule{0.5pt} \\[0.4cm] 
    \Large \bfseries #1 \\[0.2cm]
    \horrule{2pt}
  }
}
\newcommand{\Author}[2]{
    \def\Auth{#1} 
    \author{
        \vskip 10pt
        {\large{#1}} \\  \vskip -2pt
        {#2}\\
    }
}
\date{}
\newtheorem{thm}{Theorem}[section]
\newtheorem{prop}[thm]{Proposition}
\newtheorem{lem}[thm]{Lemma}
\newtheorem{cor}[thm]{Corollary}

\newtheorem{defn}[thm]{Definition}
\newtheorem{rem}[thm]{Remark}
\newcommand{\pf}{{\bf Proof. \ }}
\renewcommand{\epsilon}{\varepsilon}
\renewcommand{\rho}{\varrho}
\renewcommand{\phi}{\varphi}

\newcommand{\NN}{\ensuremath{\mathbb{N}}\xspace}
\newcommand{\ZZ}{\ensuremath{\mathbb{Z}}\xspace}
\newcommand{\QQ}{\ensuremath{\mathbb{Q}}\xspace}
\newcommand{\RR}{\ensuremath{\mathbb{R}}\xspace}
\newcommand{\CC}{\ensuremath{\mathbb{C}}\xspace}
\newcommand{\FF}{\ensuremath{\mathbb{F}}\xspace}
\newcommand{\TT}{\ensuremath{\mathbb{T}}\xspace}
\newcommand{\RP}{\ensuremath{\mathbb{RP}}\xspace}
\newcommand{\Sp}{\mathcal{S}p}
\newcommand{\TC}{\operatorname{TC}}
\newcommand{\thh}{\operatorname{THH}}
\newcommand{\bw}{B^{cy}(\Pi, \overline{\omega})}
\newcommand{\tp}{\operatorname{TP}}
\newcommand{\hofib}{\operatorname{hofib}}
\newcommand{\coker}{\operatorname{coker}}
\Title{On The Algebraic $K$-Theory of Double Points}
\Author{Noah Riggenbach}{riggenbachn@gmail.com}
\begin{document}
\maketitle

\begin{abstract}
    In this paper, we use trace methods to study the algebraic $K$-theory of rings of the form $R[x_1,\ldots, x_d]/(x_1,\ldots, x_d)^2$. We compute the relative $p$-adic $K$ groups for $R$ a perfectoid ring. In particular, we get the integral $K$ groups when $R$ is a finite field, and the integral relative $K$ groups $K_*(R[x_1,\ldots, x_d]/(x_1,\ldots, x_d)^2, (x_1,\ldots, x_d))$ when $R$ is a perfect $\FF_p$-algebra. We conclude the paper with some other notable computations, including some rings which are not quite of the above form.
\end{abstract}

\tableofcontents
\section{Introduction}
Recent years have seen a remarkable increase in interest in algebraic $K$-theory. The work of \cite{BGT} and \cite{Barwick} have described a universal property of $K$-theory. On the more calculational side, the work of Voevodsky and many others has lead to the solution to the Quillen-Lichtenbaum conjecture. We are able to say more than ever before what algebraic $K$-theory looks like.

That being said, many aspects of $K$-theory remain mysterious. A key example of this is what $K$-theory looks like at singular schemes. Many of the helpful calculational tools, such as the Quillen-Lichtenbaum conjecture and other tools coming from motivic homotopy theory, require the scheme to be regular. Once $\mathbb{A}^1$-homotopy invariance fails, we lose many of these theorems. Despite this, there have been many successful computations of the $K$-theory of singular schemes using trace methods, particularly in characteristic $p>0$ or rationally, e.g. \cite{HM_cyclic_polytopes}, \cite{Hesselholt_Madsen}, \cite{Hesselholt_Nikolaus}, \cite{LM_Extensions_by_direct_sums}, \cite{Speirs_coordinate_axes}, and \cite{Speirs_truncated_polynomials}. 

Define for a given ring $R$\[A_d:= R[x_1,\ldots, x_d]/(x_1,\ldots, x_d)^2\]and \[\mathfrak{m}:= (x_1,\ldots, x_d)\subseteq A_d.\] Then the goal of the present paper is to compute $K_*(A_d, \mathfrak{m})^\wedge_p$. We manage to do this for $R$ a perfectoid ring. In order to state the main theorem, we will need some terminology from \cite{Speirs_coordinate_axes}. This is the first of many connections to \cite{Speirs_coordinate_axes} and \cite{Speirs_truncated_polynomials}, which we explore in detail in Section~\ref{sec: Speirs}.

We take from \cite{Speirs_coordinate_axes} two functions, $t_{ev}$ and $t_{od}$. The function $t_{ev} = t_{ev}( p , r, m')$ is the unique positive integer, if it exists, such that $m'p^{t_{ev}-1}\leq 2r<m'p^{t_{ev}}$. If no such integer exists, $t_{ev}=0$. Similarly, $t_{od}=t_{od}(p,r,m')$ is, if it exists, the unique positive such that $m'p^{t_{od}-1}\leq 2r+1<m'p^{t_{od}}$, and is zero if no such integer exists. We also define the sets $\omega_{s,d}$ as the set of cyclic words, as defined in definition~\ref{defn: word}, of d letters, length $s$, and period exactly $s$. In other words, the elements of $\omega_{s,d}$ are strings of $s$ letters from an alphabet of $d$ letters such that the $C_s$ action of rotating the word acts freely. In her thesis, Rudman \cite{Rudman} has calculated the cardinality of this set as  \[|\omega_{s,d}| = \frac{\sum_{u\mid s}\mu(s/u)d^u}{s}\]where $\mu(-)$ is the M\"obius function. Finally, let $J_p=\{m'\in \ZZ_+| (m',p)=1\}$.

Our main result is now the following.
\begin{thm}\label{thm: main result}
Let $p$ be an odd prime and let $R$ be a perfectoid ring. Then, using the language above, there are isomorphisms
\[
\pi_*\left(K(A_d,\mathfrak{m})^\wedge_p\right)\cong \begin{cases}
\prod\limits_{m'\in 2J_p}\ \prod\limits_{\substack{s\mid m'p^{t_{ev}-1}\\s\textrm{ even}}}\ \prod\limits_{\omega_{s,d}}W_{t_{ev}-v_p(s)}(R)  &   \textrm{if }*=2r\\
\prod\limits_{m'\in J_p\setminus 2J_p}\ \prod\limits_{s\mid m'p^{t_{od}-1}}\ \prod\limits_{\omega_{s,d}}W_{t_{od}-v_p(s)}(R)   &   \textrm{if }*=2r+1
\end{cases}
\]
where $W_n(R)$ are the truncated $p$-typical Witt vectors of $R$ and $v_p$ is the $p$-adic valuation.

For $p=2$ and $R$ perfectoid with respect to $p$, we have isomorphisms
\[\pi_*\left(K(A_d,\mathfrak{m})^\wedge_2\right)\cong \begin{cases}
\prod\limits_{m'\in 2\ZZ}\prod\limits_{\substack{s\mid m\\ s \textrm{ even}}}\prod\limits_{\omega_{s,d}}W_{t_{ev}-v_{2}(s)}(R)  &   \textrm{if }*=2r\\
\prod\limits_{m'\in \ZZ\setminus 2\ZZ}\prod\limits_{s\mid m'}\ \prod\limits_{\nu=0}^{t_{ev}-1}\prod\limits_{\omega_{s,d}}R &   \textrm{if }*=2r+1
\end{cases}.\]
\end{thm}

\begin{rem}
In Theorem~\ref{thm: main result}, the products are finite. This is because for $m'$ large enough $t_{ev}=t_{od}=0$.
\end{rem}

For the reader unfamiliar with perfectoid rings, the relevant theory is reviewed in Section~\ref{sec:perfectoid rings}. They key example to keep in mind is when $R$ is a perfect $\FF_p$-algebra, i.e., it is an $\FF_p$-algebra where the Frobenius is an isomorphism. The main results cited and proven in Section~\ref{sec:perfectoid rings} are lifts of results known for perfect $\FF_p$-algebras to the setting of perfectoid rings.

This work originally started as a revisit of \cite{LM_Extensions_by_direct_sums}, with the goal of reproducing the calculation using the new methods introduced by Nikolaus and Scholze in \cite{Nikolaus_Scholze}. We discuss how Theorem~\ref{thm: main result} relates to \cite{LM_Extensions_by_direct_sums} briefly in Section~\ref{sec: perfect algebra}, where we extend their result from perfect fields of characteristic $p>0$ to perfect $\FF_p$-algebras. We have also managed to generalize some of the work of Speirs, as discussed in Section~\ref{sec: Speirs}.

\subsection{Outline} We prove Theorem~\ref{thm: main result} using trace methods. Section~\ref{sec: K to TC} is the reduction to a $\TC$ calculation using trace methods. Section~\ref{sec: THH} then computes the topological Hochschield Homology as a  $\TT$-equivariant spectrum. It turns out that this spectrum is comprised of spectra with induced action at odd primes, and at $p=2$ is only slightly different, and so Section~\ref{sec: TC- and TP} is a computation of topological negative cyclic and periodic homology in terms of the homotopy fixed points and Tate construction for the cyclic subgroups of $\TT$. Section~\ref{sec: tc} combines this along with technical results from Section~\ref{sec: triv spheres} to get Theorem~\ref{thm: main result} using \cite{Nikolaus_Scholze}. Finally, Section~\ref{sec: examples} contains some consequences of Theorem~\ref{thm: main result}. For a similar strategy, see \cite{Hesselholt_Nikolaus}, and for a spectral sequence based approach to similar problems, see \cite{Speirs_coordinate_axes} or \cite{Speirs_truncated_polynomials}.

\subsection{Acknowledgements} I would like to thank my advisor Michael Mandell for giving me useful feedback on this paper, and for being a helpful and supportive advisor. I would also like to extend a special thanks to Ayelet Lindenstrauss and Emily Rudman, who both helped me work through many parts of this paper, pointed out multiple mistakes I made, and made the writing significantly better. I would like to thank Sanjana Agarwal and Martin Speirs, for patiently helping me understand perfectoid rings, and for many other helpful conversations. I would also like to thank Martin Speirs for many helpful comments on an earlier version of this paper, including pointing out a crucial mistake which Section~\ref{sec: triv spheres} is added to address. Finally, I would like to thank Dylan Spence for giving an algebraic geometer's opinion of Section~\ref{sec:perfectoid rings}. 

While working on this paper I was supported by a Hazel King Thompson Scholarship from the Mathematics Department at Indiana University.

\section{\texorpdfstring{Reduction to $\TC$ calculation}{Reduction to TC calculation}}\label{sec: K to TC}

The ring $A_d$ is particularly well suited for this computation, since on it we may apply the Dundas Goodwillie McCarthy Theorem \parencite[Theorem 7.0.0.2]{Dundas_Goodwillie_McCarthy} (or for a more modern approach and statement, \parencite[Theorem 1.1.1]{Raskin}) to get a pullback square

\begin{center}
\begin{tikzcd}
K(A_d) \arrow[r] \arrow[d]                  & \TC(A_d) \arrow[d]                   \\
K(R) \arrow[r] \arrow[u, dotted, bend left] & \TC(R) \arrow[u, dotted, bend right]
\end{tikzcd}
\end{center}
where the vertical maps are the quotient by the (nilpotent) ideal $\mathfrak{m}$. As noted in the diagram, the inclusion $R\to A_d$ gives a splitting of the above maps, and hence we get an isomorphism 
\begin{equation}
K(A_d)\cong K(R)\vee K(A_d, \mathfrak{m})\cong K(R)\vee \TC(A_d, \mathfrak{m})\label{eqn:k-split}
\end{equation}
where the last equivalence comes from the above pullback square.

The celebrated computation of Quillen \cite{Quillen} gives the homotopy groups of the first summand in the case of $R$ a finite field. In the more general case of $R$ a perfectoid ring, we will usually not be able to say more in the integral setting.

Our main theorem is then equivalent to the following:

\begin{thm}\label{thm: TC main thm}
For $p$ an odd prime and $R$ a perfectoid ring, there are isomorphisms \[\pi_{2r}\left(\TC(A_d, \mathfrak{m})^\wedge_p\right)\cong \prod_{m'\in 2J_p}\ \prod_{\substack{s\mid m'p^{t_{ev}(p, r, m')-1}\\s\textrm{ even}}}\ \prod_{\omega_{s,d}}W_{t_{ev}(p, r, m')-v_p(s)}(R)\]
and
\[\pi_{2r+1}\left(\TC(A_d, \mathfrak{m})^\wedge_p\right)\cong \prod_{m'\in J_p\setminus 2J_p}\ \prod_{s\mid m'p^{t_{od}(p, r, m')-1}}\ \prod_{\omega_{s,d}}W_{t_{od}(p, r, m')-v_p(s)}(R)\]
where $v_p(-)$ is the $p$-adic valuation.

For $p=2$ and $R$ a perfectoid ring with respect to $2$, there are isomorphisms
\[\pi_{2r}\left(\TC(A_d,\mathfrak{m})^\wedge_2\right)\cong \prod_{m'\in 2\ZZ}\ \prod_{\substack{s\mid m'2^{t_{ev}-1}\\ s\textrm{ even}}}\prod_{\omega_{s,d}}W_{t_{ev}-v_2(s)}(R)\]
and
\[
\pi_{2r+1}\left(\TC(A_d,\mathfrak{m})^\wedge_2\right)\cong \prod_{m'\in \ZZ\setminus 2\ZZ}\prod_{s\mid m'}\prod_{\nu=0}^{t_{ev}-1}\prod_{\omega_{s,d}}R.
\]
\end{thm}

\begin{rem}
If we are willing to work rationally, Goodwillie's theorem \cite{Goodwillie} lets us replace~Equation (\ref{eqn:k-split}) with $$K(A_d)_\QQ\cong K(R)_\QQ \vee \Sigma \operatorname{HC}(A_d\otimes \QQ, \mathfrak{m}\otimes \QQ).$$ The second summand was studied in Rudman \cite{Rudman} in the case of $R=\ZZ$, which computes $\operatorname{HC}(A_d)$ for more general $R$. In addition,  Thomason \cite{Thomason}  showed that $K(-)_\QQ$ satisfies \'etale descent for schemes under mild hypotheses, and is more amenable to computation then ordinary algebraic $K$-theory. 
\end{rem}

For the remainder of this section we will work in the $p$-complete setting. In this setting we see one of the first benefits of the recent advances in perfectoid rings and homotopy theory.

\begin{thm}[\cite{Clausen_Mathew_Morrow}, Theorem B]
Let $R$ be a ring henselian along $(p)$ and such that $R/p$ has finite Krull dimension. Let $d = \sup_{x\in \operatorname{Spec}(R/p)}\log_p[k(x):k(x)^p]$, where $k(x)$ denotes the residue field at $x$. Then the map $K(R)/p^i\to \TC(R)/p^i$ is an equivalence in degrees $\geq \operatorname{max}(d,1)$ for all $i\geq 1$.
\end{thm}

As noted in \cite{Clausen_Mathew_Morrow}, their proof of this result specializes nicely to the case of $R$ semiperfectoid. For our purposes we will only need the following corollary.

\begin{cor}[\cite{Clausen_Mathew_Morrow}]
Let $R$ be a perfectoid ring. Then the map $K(R)^\wedge_p\to \TC(R)^\wedge_p$ exhibits the former as the connective cover of the latter.
\end{cor}

Hence going from a calculation of the relative $K$-theory $K(A_d,\mathfrak{m})^\wedge_p$ to a computation of $K(A_d)^\wedge_p$ is a matter of computing $\TC(R)^\wedge_p$ for these rings. In principal, this is an easier computation. In particular, the recent work of Bhatt and Scholze in \cite{BS} introduces a promising calculational tool for $R$ quasiregular semiperfectoid.

To finish this section, we note that we do not lose any information in $p$-completion when $R$ is a perfect $\FF_p$-algebra.
\begin{prop}\label{prop: perfect p complete}
Let $R$ be a perfect $\FF_p$-algebra. Then all of the homotopy groups of $\TC(A_d, \mathfrak{m})$ are $p$-power torsion. In particular, $\TC(A_d, \mathfrak{m})$ is $p$-complete. 
\end{prop}
\pf Since $p=0$ in $A_d$, $\mathfrak{m}$ is $p$-power torsion. Applying \parencite[Theorem D]{Land_Tamme} then gives the result.\qed 

\section{\texorpdfstring{Computation of $\thh$ as a $\TT$-spectrum}{Computation of THH as a T-spectrum}}\label{sec: THH}

We begin by noting that $A_d = R\wedge \Pi$, where $\Pi = \{0, 1, x_1,\ldots, x_d\}$ is the pointed monoid with all products $x_ix_j = 0$. Since $\thh$ is symmetric monoidal \parencite[Section IV.2]{Nikolaus_Scholze}, we then get that $$\thh(A_d)\simeq \thh(R)\wedge \thh(\Sigma^\infty \Pi).$$ The latter term is equivalent to $\Sigma^\infty B^{cy}(\Pi)$. We refer the reader to \parencite[Sections 3.1 and 3.2]{Speirs_coordinate_axes} for a comprehensive review of the relevant details about this decomposition. In particular, we have the following lemma.

\begin{lem}[\cite{Hesselholt_Madsen}, Theorem 7.1; \cite{Nikolaus_Scholze}, Section IV.2]\label{lem:smash}
Let $R$ be a ring, $\Pi$ a pointed monoid, and $R[\Pi]$ the pointed monoid algebra. Then there is a natural $\TT$-equivariant equivalence $$\thh(R)\wedge B^{cy}(\Pi) \xrightarrow{\sim} \thh(R[\Pi])$$ where the $\TT$ action on the left is the diagonal action. Under this equivalence, the Frobenius map is induced by the Frobenius on $\thh(R)$ smashed with the unstable Frobenius on $B^{cy}(\Pi)$.
\end{lem}

By naturality, we see that the map $\thh(R)\to \thh(R[\Pi])$ is given by the composition $$\thh(R)\cong \thh(R)\wedge S^0 \to \thh(R)\wedge B^{cy}(\Pi).$$ Consequently, the cofiber of this map is given by $\thh(R)\wedge \operatorname{cofib}(S^0\to B^{cy}(\Pi))\cong  \thh(R)\wedge \widetilde{B^{cy}(\Pi)}$, where $\widetilde{B^{cy}(\Pi)}$ is the cofiber of the map $S^0\to B^{cy}(\Pi)$ in spaces. 

\subsection{\texorpdfstring{Weight and Cyclic Word Decomposition of $B^{cy}(\Pi)$}{Weight and Cyclic Word Decomposition of the Cyclic Bar Construction}}

As described in \parencite[Section 3.2]{Speirs_coordinate_axes}, $B^{cy}(\Pi)$ can be decomposed as a pointed $\TT$-space into a wedge of simpler spaces. In order to do this decomposition, we need to define the notions of word, word length, word period, and cyclic word.

\begin{defn}\label{defn: word}
Consider the set $S= \{x_{1},\ldots, x_{d}\}$. We define a word $\omega$ of length $m\geq 0$ to be a mapping  $\omega: \{1,2,\ldots,m\} \rightarrow S.$ The cyclic group $C_m$ acts on the words of length $m$; the orbit of a word $\omega$ is denoted $\overline{\omega}$ and is called a cyclic word. The cardinality of $\overline{\omega}$ is the period of $\omega$ which will we denote as $s$. For $m>0$, note for any word of length $m$, $s$ must divide $m$. For $m=0$, $s=1$.
\end{defn}

The n-simplicies $B^{cy}_n(\Pi)$ of $B^{cy}(\Pi)$ correspond to the $(n+1)$-tuples $(a_{0},\ldots , a_{n})$ where $a_{i} \in \{1, x_{1},\ldots, x_{d}\}$ for $i \in \{0, 1, \ldots, n\}$. Each $(n+1)$-tuple corresponds to exactly one cyclic word $\overline{\omega}$ by ignoring any $1$'s in the sequence. The face maps $d_{0}, \ldots, d_{n-1}$ and degeneracy maps $s_{0}, \ldots, s_{n}$ all will send the $n$-simplex corresponding to $(a_{0}, \ldots , a_{n})$ either to the basepoint 0 or to a simplex corresponding the same reduced sequence. The $n^{th}$ face map $d_{n}$ and the cyclic operator $t_{n}$ will both send  $(a_{0}, \ldots , a_{n})$ either to $0$ or to a simplex whose corresponding reduced cyclic word is $\overline{\omega}$. Since all the structure maps preserve $\overline{\omega}$ we get a decomposition

\begin{equation}\label{Bsplit}  B^{cy}(\Pi) = \bigvee_{\textrm{cyclic words }\overline{\omega} } B^{cy}(\Pi, \overline{\omega})
\end{equation}

$\\*$where $B^{cy}(\Pi, \overline{\omega})$ is the subspace of $B^{cy}(\Pi)$ whose cells all correspond to the same $\overline{\omega}$ or to $0$ for any $\overline{\omega}$ that is a cyclic word in the letters of $S$. Note that when $m=0$, we are exactly getting the simplicies with only zeros and ones. Hence the summand corresponding to $m=0$, i.e. the empty cyclic word, is homeomorphic to $S^0$.

\subsection{\texorpdfstring{Equivariant Homotopy Type of $\bw$}{Equivariant Homotopy Type of the Cyclic Word Decomposition}}

\begin{lem}\label{OurMartinLem3.3}
Let $\overline{\omega}$ be a cyclic word of cycle length $s\geq 2$, with letters in the set $S=\{x_{1}, x_{2},\hdots, x_{d}\}$ and with length $m=s\cdot i$. A choice of $\omega$ with cyclic word $\overline{\omega}$ determines a $\TT$-equivariant homeomorphism $$S^{\RR[C_{m}]-1}\wedge_{C_{i}}\TT_{+} \xrightarrow{\sim} B^{cy}(\Pi, \overline{\omega})$$ where $\RR[C_{m}]-1$ is the reduced regular representation of $C_{m}.$
\end{lem}

\pf
This follows almost directly from the proof of \parencite[Lemma 9]{Speirs_coordinate_axes} noting the difference in pointed monoids. In the paper \parencite{Speirs_coordinate_axes}, 
$\Pi^{d} = \{0,1,x_{1}, x_{1}^{2}, \ldots, x_{2}, x_{2}^{2}, \ldots, x_{d}, x_{d}^{2}, \ldots \}$ is the multiplicative monoid with basepoint 0 and multiplication $x_{i}x_{j} = 0$ for $i\neq j$.  He defines a word of length $m$ with  no cyclic repetitions \parencite[Definition 8]{Speirs_coordinate_axes} to be a word $\omega = w_{1}w_{2}\ldots w_{m}$ such that $w_{i}\neq w_{i+1}$ for $i=0,1,\ldots,m-1$ and $w_{m} \neq w_{1}$. 

The proof for \parencite[Lemma 9]{Speirs_coordinate_axes} for cyclic words $\overline{\omega}$ of length $m$ with no cyclic repetitions relies on the following: all of the faces of any $m$-simplex in $B^{cy}(\Pi,\overline{\omega})$ with corresponding reduced word in $\overline{\omega}$ are identified to the basepoint $0$.   In the pointed monoid used for the coordinate axes paper, this is because attention was restricted to words with no cyclic repeats. However, for our pointed monoid, $x_{i}x_{j}= 0$ for any $i,j$.  Therefore, all faces of any $m$-cell in $B^{cy}(\Pi,\overline{\omega})$ with corresponding reduced word in $\overline{\omega}$ will be $0$ regardless of whether the word $\omega$ has cyclic repetitions or not. The rest of the proof follows  directly from the proof of \parencite[Lemma 9]{Speirs_coordinate_axes}.
\qed

\begin{lem}
Let $\overline{\omega}$ be a cyclic word of cycle length $s\geq 1$, with letters in the set $S=\{x_{1}, x_{2},\hdots, x_{d}\}$ and with length $m=s\cdot i$.

\begin{enumerate}
    \item For $s$ even, there is a $\TT$-equivariant homeomorphism
    $$
    \Sigma B^{cy}(\Pi,\overline{\omega})\cong S^{\lambda_{m/2}}\wedge \left(\TT/C_{i}\right)_{+}.
    $$
    
    \item For $s$ and $i$ both odd, there is a $\TT$-equivariant homeomorphism
    $$
    B^{cy}(\Pi,\overline{\omega})\cong S^{\lambda_{(m-1)/2}}\wedge \left(\TT/C_{i}\right)_{+}.
    $$
    \item For $s$ odd and $i$ even, there is a $\TT$-equivariant homeomorphism
    $$
    B^{cy}(\Pi,\overline{\omega})\cong S^{\lambda_{(m-2)/2}}\wedge \RP^{2}(i).
    $$
\end{enumerate}

\end{lem}

\smallskip
Here $\mathbb{R}P^2(i)$ is the cofiber of the map $\left(\mathbb{T}/C_{\frac{i}{2}}\right)_+\to \left(\mathbb{T}/C_i\right)_+$ in $\TT$-spaces, and $\lambda_n\cong \CC(1)\oplus \CC(2)\oplus\cdots\oplus\CC(n)$, where $\CC(j)$ is the one-dimensional complex representation of $\TT$ having $z\in \TT$ act by multiplication by $z^j$.

\pf
See Lemma \ref{OurMartinLem3.3} and \parencite[Lemma 11]{Speirs_coordinate_axes}.\qed

\section{Trivializing the representation spheres}\label{sec: triv spheres}

In light of the previous section, we get the following decomposition of $\thh(A_d)$:
\begin{equation}
\thh(A_d) \simeq \thh(R)\vee\bigvee_{m\in \ZZ_+}
\begin{bmatrix}
\bigvee_{\substack{s\mid m\\ s\textrm{ even}}}\left( \bigvee_{\overline{\omega}\in \omega_{s,d}}\Sigma^{-1}\thh(R)\wedge S^{\lambda_{m/2}}\wedge \left(\TT/C_{\frac{m}{s}}\right)_+\right)\hfill\\
\vee \bigvee_{\substack{s\mid m\\ m \textrm{ odd}}}\left(\bigvee_{\overline{\omega}\in \omega_{s,d}} \thh(R)\wedge S^{\lambda_{m-1/2}}\wedge \left(\TT/C_{\frac{m}{s}}\right)_+\right)\hfill\\
\vee \bigvee_{\substack{s\mid m\\ s\neq m \textrm{ mod } 2}}\left(\bigvee_{\overline{\omega}\in \omega_{s,d}} \thh(R)\wedge S^{\lambda_{m-2/2}}\wedge \RP(\frac{m}{s})\right)\hfill
\end{bmatrix}\label{eqn: Wedge decomp of THH}
\end{equation}
where $\omega_{s,d}$ is the set of all cyclic words with cycle length actually equal to $s$ (and not a proper divisor of $s$) in $d$ letters.
In addition, the map $\thh(R)\to \thh(A_d)$ giving the splitting of $\TC$ is the inclusion as the first summand.

Note also that each summand is of the form $\thh(R)\wedge S^\lambda \wedge X$ for some representation $\lambda$ and space $X$. The $\TT$-action on this can be quite complicated if we allow arbitrary $\mathbb{E}_1$-ring spectra $R$. We make the following definition in order simplify these spaces.

\begin{defn}\label{def: T-sus inv}
Let $X$ be any $\TT$-equivariant spectrum. We say that $X$ is $\TT$-suspension invariant if for any $\TT$-representation sphere $S^V$, there is a Borel equivalence \[X\wedge S^V\simeq \Sigma^{\dim(V)}X.\]
\end{defn}

We will now spend some time showing that topological hochschield homology of perfectoid rings satisfy this condition. In fact, we will show something stronger:

\begin{prop}\label{prop: even spectra trivial action}
Let $X$ be any element of $E\textrm{-mod}^{B\TT}$, $E$ an $\mathbb{E}_\infty$-ring spectrum such that $\pi_{2r+1}(E)=0$ for all $r\in \ZZ$. Then $X$ is $\TT$-suspension invariant.
\end{prop}

We will see in Corollary~\ref{cor: T-sus independance of even dodads} that this Proposition applies in our situation, and so we may consider only trivial representation spheres. Before we prove this result, we first need the following lemma.

\begin{lem}
Let $E$ be an $\mathbb{E}_\infty$-ring spectrum with $\pi_{2r+1}(E)=0$ for all $r\in \ZZ$. Let $X$ be any element of $E\textrm{-mod}^{B\TT}$ such that $F(X)\simeq \Sigma^n E$ for some $n\in \ZZ$, $F:E\textrm{-mod}^{B\TT}\to E\textrm{-mod}$ the forgetful functor. Then there is a Borel equivalence $X\simeq (\Sigma^n E)^{triv}$, here $(-)^{triv}$ is the trivial action functor. 
\end{lem}

\pf
Consider $X$ as an $(\infty, 1)$-functor $X:B\TT\to E\textrm{-mod}$. It is enough to show that $X$ is contractable as such a map of pointed simplicial sets (pointing $E\textrm{-mod}$ by $X(*)$), since then the identification map $\Sigma^n E \simeq X$ will be equivariant.

This map then factors through a map $B\TT\to Bhaut_{E}(X(*))$, since the sources is a groupoid. Here $haut_E(-)$ is the connected component of the space of $E$-module maps which contain a weak equivalence. Furthermore, since both spaces are Kan complexes, it is then enough to show after geometric realization that this map is contractable. Under the identification of $X(*)\simeq \Sigma^n E$, this is then equivalent to a map $B\TT\to Bhaut_E(\Sigma^n E)$. By (de)suspending by trivial representation spheres, we may assume without loss of generality that $n=0$. It remains to classify that maps $B\TT\to Bhaut_E(E)$.

To this end, note that $haut_E(E)\subseteq E\textrm{-mod}(E,E)\simeq \Omega^{\infty}E$\footnote{Warning: This identification is only as a space, and not as an $\mathbb{E}_\infty$-space. We go on to take $B$ of the $\mathbb{A}_\infty$ structure given by composition and not the lax-monoidal struture of $E\textrm{-mod}(E,-)$, so this distinction matters.}, and since $haut_E(E)$ is a subspace of connected components, $\pi_i(haut_E(E))\cong \pi_i(\Omega^\infty E)\cong \pi_i(E)$ for $i\geq 1$. In particular, $haut_E(E)$ has homotopy groups concentrated in even degree.

Note also that $E\textrm{-mod}$ is a symmetric monoidal $(\infty,1)$-category, and so the monoidal unit will have an $\mathbb{E}_\infty$ structure on its homotopy automorphism space, and the standard $\mathbb{A}_\infty$ structure will be the one induced by restricting this $\mathbb{E}_\infty$ structure. We then have a group-like $\mathbb{E}_\infty$-space $haut_E(E)$, and so $\pi_*(Bhaut_E(E))\cong \pi_{*-1}(haut_E(E))$ is concentrated in odd degree.

In addition, $Bhaut_E(E)$ is also and $\mathbb{E}_\infty$-space, and so homtopy classes of maps into it can be computed using the Atiyah Hirzebruch spectral sequence:
\[E_2^{p,q}=\widetilde{H}^p(B\TT; \pi_{-q}(Bhaut_E(E)))\implies [B\TT, Bhaut_E(E)]_{-p-q}.\]
The $E_2$ page is concentrated in odd total degree, and there is no room for differentials so the spectral sequence converges. Thus there are no non-contractable maps of even degree, and in particular no non-contractable maps $B\TT\to Bhaut_E(E)$.\qed

We are now ready to prove Proposition~\ref{prop: even spectra trivial action}.

\textbf{Proof of Proposition~\ref{prop: even spectra trivial action}. \ } Suppose $X$ is any element of $E\textrm{-mod}^{B\TT}$, and let $S^V$ be any $\TT$-representation sphere. Then we have the following chain of Borel equivalences:
\[X\wedge S^V\simeq \left(X\wedge_E E\right)\wedge S^V\simeq X\wedge_E(E\wedge S^V).\] By the lemma above, there is a Borel equivalence $E\wedge S^V\simeq \Sigma^{\dim(V)}E$ in $E\textrm{-mod}^{B\TT}$, where the spectrum on the right has trivial $\TT$-action. Hence \[X\wedge_E\left(E\wedge S^V\right)\simeq X\wedge_E\left(E\wedge S^{\dim(V)}\right)\simeq X\wedge S^{\dim(V)}\] in $E\textrm{-mod}^{B\TT}$, as desired.\qed

In particular, we get the following useful corollary, which in particular applies to $\thh(R)^\wedge_p$ for $R$ quasiregular semiperfectoid (and hence also perfectoid).

\begin{cor}\label{cor: T-sus independance of even dodads}
Let $X$ be any $\TT$-equivariant $\mathbb{E}_\infty$-ring spectrum which is concentrated in even degrees. Then $X$ is $\TT$-suspension invariant.
\end{cor}
\pf
Since $X$ is concentrated in even degrees, the homotopy fixed point spectral sequence converging conditionally to the homotopy groups of $X^{h\TT}$ has $E_2$ page concentrated in even bidegrees. Hence the spectral sequence collapses at the $E_2$ page, converges, and has no nonzero entries in odd total degree. Consequentially $X$ and $E=X^{h\TT}$ satisfy the hypotheses for Proposition~\ref{prop: even spectra trivial action}.\qed

\section{Topological negative cyclic and periodic homology}\label{sec: TC- and TP}

For $R$ a perfectoid ring, Corollary~\ref{cor: T-sus independance of even dodads} applies and $\thh(R)^\wedge_p$ is $\TT$-suspension invariant (see Definition~\ref{def: T-sus inv}.). Consequentially, we may take the trivial representation spheres in Equation~\ref{eqn: Wedge decomp of THH}. Note that since the representations $\lambda_i$ have complex dimension $i$ (so real dimension $2i$), these trivial representation spheres are $S^{n}$, where $n$ is $m-1$ or $m-2$, depending on the case. In particular, if we write $\thh(A_d)= \bigvee_{m\in \NN} \thh(m)$ where $\thh(m)$ is the wedge of all summands corresponding to cyclic words of length $m$, then the connectivity of $\thh(m)$ is at least $m-2$ and so
\begin{equation}
\thh(A_d)\simeq \bigvee_{m\in \NN}\thh(m)\simeq \prod_{m\in \NN}\thh(m).\label{eqn: product decomp of THH}
\end{equation}

We may further decompose $\thh(A_d)$ into a product by noting that in Equation~(\ref{eqn: Wedge decomp of THH}), each of the inner wedge sums are finite wedges so they are also equivalent to the respective products. Therefore, to compute $\TC^{-}(A_d)$, and the relative $\TC^{-}$, it is enough to compute $(-)^{h\TT}$ on each of the the summands in Equation ~(\ref{eqn: Wedge decomp of THH}), and $(-)^{h\TT}$ for the summands where $m\geq 1$, respectively. 

We may make a similar observation for the homotopy orbits as well. Since Equation~(\ref{eqn: Wedge decomp of THH}) is given in terms of iterated wedge sums, $(-)_{h\TT}$ can be computed term by term. In addition, homotopy orbits only increase connectivity, so in the decomposition $\thh(A_d)_{h\TT}\simeq \bigvee_{m\in \NN}\thh(m)_{h\TT}$, each term $\thh(m)_{h\TT}$ is at least $(m-2)$-connected. Consequently, we may express $(\thh(A_d))_{h\TT}$ as a product in the exact same way we did for $\TC^{-}(A_d)$, and the canonical map will respect this decomposition. In particular, we get the same decomposition of the topological periodic homology. 

Note that each summand in Equation~(\ref{eqn: Wedge decomp of THH}) with $s=m\mod 2$ is of the form $X\wedge \left(\TT /C_n\right)_+$. We use the following proposition to simplify the computation.

\begin{prop}[\cite{Hesselholt_Nikolaus}, Proposition 3]\label{prop: Induced/coinduced action and fixed points}
Let $G$ be a compact Lie group. Let $H\subseteq G$ be a closed subgroup, let $\lambda=T_H(G/H)$ be the tangent space at $H=eH$ with the adjoint left $H$-action, and let $S^\lambda$ be the one-point compactification of $\lambda$. For every spectrum with $G$-action $X$, there are canonical natural equivalences \[\left(X\wedge \left(G/H\right)_+\right)^{hG}\simeq \left(X\wedge S^\lambda\right)^{hH},\] \[\left(X\wedge \left(G/H\right)_+\right)^{tG}\simeq \left(X\wedge S^\lambda\right)^{tH}.\]
\end{prop}
In particular, $\left(X\wedge (\TT/C_n)_+\right)^{h\TT}\simeq (\Sigma X) ^{hC_n}$, and similarly for the Tate construction. Furthermore, the Tate construction is initial among functors under the homotopy fixed points vanishing on compact objects by \parencite[Lemma I.1.4(ii)]{Nikolaus_Scholze}, and so the canonical map must be sent to the canonical map under this equivalence. 

\begin{rem}
The proof of Proposition \ref{prop: Induced/coinduced action and fixed points} in \cite{Hesselholt_Nikolaus} is from the point of view of infinity categories, and uses constructions from \cite{Nikolaus_Scholze}. While the proof using infinity categories has simplified things considerably,  \ref{prop: Induced/coinduced action and fixed points} was known before.  It can be obtained from the Wirthm\"uller isomorphism, which is proven in \cite{May} in a similar fashion as above, but in terms of a six functor formalism of model categories from \cite{Fausk_Hu_May}. 
\end{rem}

\subsection{\texorpdfstring{Recollection of the calculation of topological Hochschild homology of perfectoid rings}{Recollection of topological hochschild homology of perfectoid rings}}\label{sec:perfectoid rings}

In order to do explicit calculations, we need to know more about $\thh(R)$. For a review on the algebraic properties of perfectoid rings, we refer the reader to \cite{BMS1}. We recall the definition of perfectoid rings here for convenience:

\begin{defn}[\cite{BMS1}, Definition 3.5]\label{defn: perfectoid}
A ring $R$ is perfectoid if there exists an element $\pi \in R$ such that \begin{enumerate}
    \item $\pi^p$ divides $p$;
    \item $R$ is $\pi$-adically complete, and is separated with respect to this topology (and therefore p-complete);
    \item The Frobenius map $\phi: R/p\to R/p$ is surjective;
    \item The kernel of the map $\theta: A_{\textrm{inf}}(R)\to R$ is principal.
\end{enumerate}
where $A_{\textrm{inf}}(R)$ is Fontaine's ring, $W\left(\lim\left( \ldots \xrightarrow{\phi}R/p\xrightarrow{\phi} R/p\right)\right)$ that is: the ($p$-typical) Witt vectors on $\lim\left( \ldots \xrightarrow{\phi}R/p\xrightarrow{\phi} R/p\right)$.
\end{defn}
Following \cite{BS}, we will refer to a choice of element $\xi\in A_{\textrm{inf}}(R)$, $(\xi)=ker(\theta)$, as an orientation of $A_{\textrm{inf}}(R)$. By abuse of notation, we will refer to such a $\xi$ as an orientation of $R$ as well. When there is no danger of confusion, we will drop the $R$ in $A_{\textrm{inf}}$.

Notice that since $A_{\textrm{inf}}$ is defined as the Witt vectors of a perfect $\FF_p$-algebra, it in particular comes with a Frobenius automorphism $\phi: A_{\textrm{inf}}\to A_{\textrm{inf}}$ lifting the Frobenius $\phi: R/p \to R/p$. The map $\theta:A_{\textrm{inf}}\to R$ is one of a family of maps $\theta_r:A_{\inf}\to W_{r}(R)$ whose construction is reviewed in \cite{BMS1}. These maps, along with the maps $\Tilde{\theta}_r = \theta_r\circ \phi^{-r}$, are characterized by a universal property.

As noted at the end of Section~\ref{sec: K to TC}, many important theorems have been proven recently about $K$-theory and $\TC$ of these rings, at least in the $p$-complete setting, such as the extension of B\"okstedt's periodicity result:

\begin{thm}[\cite{BMS2}, Theorem 6.1]
For a perfectoid ring $R$,  $\pi_*\left(\thh(R)^\wedge_p\right)\cong R[u]$ is a polynomial ring, where $u\in \pi_2\left(\thh(R)^\wedge_p\right)\cong \pi_2 \left(\operatorname{HH}(R)^\wedge_p\right) \cong ker(\theta)/ker(\theta)^2$ can be chosen to be any generator of $ker(\theta)/ker(\theta)^2$. 
\end{thm}

It then follows that all the differentials in both the homotopy $\TT$-fixed point and Tate $\TT$-fixed point spectral sequences must be zero, since the nonzero elements are all concentrated in even dimensions. Both spectral sequences converge strongly, and the computation for $R=\FF_p$ by \cite{Nikolaus_Scholze}, \cite{BMS2} extend to any perfectoid ring $R$. With some additional work, \cite{BMS2} show the following.

\begin{prop}[\cite{BMS2}, Proposition 6.2]\label{prop: hmtpy groups of T fixed points}
The commutative square
\[
\begin{tikzcd}
\TC^{-}(R)^\wedge_p \arrow[d, "can"] \arrow[r, "\phi_p^{h\TT}"] & \operatorname{TP}(R)^\wedge_p \arrow[d, "can"] \\
\thh(R)^\wedge_p \arrow[r, "\phi_p"]                              & \thh(R)^{tC_p}                                         
\end{tikzcd}\]
gives, upon taking homotopy groups, the square
\[
\begin{tikzcd}
{A_{\textrm{inf}}[u,v]/(uv-\xi)} \arrow[rr, "\phi\textrm{-linear}"'] \arrow[rr, "{u\mapsto \sigma, v\mapsto \phi(\xi)\sigma^{-1}}"] \arrow[d, "{u\mapsto u, v\mapsto 0}"] &  & {A_{\textrm{inf}}[\sigma, \sigma^{-1}]} \arrow[d, "\sigma\mapsto \sigma"] \\
{R[u]} \arrow[rr, "u\mapsto \sigma"]                                                                                                  &  & {R[\sigma, \sigma^{-1}]}                                                 
\end{tikzcd}\]
where $\xi$ is an orientation of $R$ and the columns also apply $\theta: A_{\textrm{inf}}(R)\to R$.
\end{prop}

In addition to understanding $\TC^{-}(R)^\wedge_p$ and $\tp(R)^\wedge_p$, we also need to understand what happens at the  finite cyclic subgroups, i.e. $\left(\thh(R)^\wedge_p\right)^{hC_{p^n}}$ and $\thh(R)^{tC_{p^n}}$.  The rest of the section will be devoted to this calculation.

Combining the calculation of the Tate construction from \parencite[Remark 6.6]{BMS2} and \parencite[Lemma 3.12]{BMS1} gives \[\pi_*\left((\thh(R)^\wedge_p)^{tC_{p^n}}\right)\cong W_n(R)[\sigma, \sigma^{-1}]
\cong A_{\textrm{inf}}(R)[\sigma, \sigma^{-1}]/(\Tilde{\xi}_n)
\cong A_{\textrm{inf}}(R)[\sigma, \sigma^{-1}]/(\Tilde{\xi}_n\sigma^{-1})\]
where $\Tilde{\xi}_n = \phi(\xi)\phi^2(\xi)\ldots\phi^{n}(\xi)$, the key step of the computation is briefly reviewed in Lemma~\ref{lem: Versch map}. We  use this to compute $\pi_*\left(\left(\thh(R)^\wedge_p\right)^{hC_{p^n}}\right)$, starting with the following proposition.

\begin{prop}\label{prop: frob and can computation}
There is a choice of isomorphism such that 
the map \[ \tau_{\geq -1}\left(\thh(R)^\wedge_p\right)^{hC_{p^n}}\xrightarrow{\phi_p}
 \tau_{\geq -1}\left(\thh(R)^{tC_p}\right)^{hC_{p^{n}}}\simeq \tau_{\geq -1}\thh(R)^{tC_{p^{n+1}}}\]
(the last equivalence coming from \parencite[Lemma II.4.1]{Nikolaus_Scholze}) gives on homotopy groups the map \[A_{\textrm{inf}}[u]/(\Tilde{\xi}_{n+1})\xrightarrow{}A_{\textrm{inf}}[\sigma]/(\Tilde{\xi}_{n+1})\]
sending $u\mapsto \sigma$.
Under the same isomorphism, the map \[\tau_{\geq -1}\left(\thh(R)^\wedge_p\right)^{hC_{p^n}}\xrightarrow{can^{hC_{p^{n-1}}}}\tau_{\geq -1}\left(\thh(R)^{tC_p}\right)^{hC_{p^{n-1}}}\simeq \tau_{\geq -1}\thh(R)^{tC_{p^{n}}}\]
gives on homotopy groups the map \[W_{n+1}(R)[u] \xrightarrow{} W_{n}(R)[\sigma]\]
which reduces $W_{n+1}(R) \to W_{n}(R)$ and sends $u\mapsto \Tilde{\theta}_r(\xi)\sigma$.
\end{prop}
\pf By Proposition~\ref{prop: hmtpy groups of T fixed points}, the Frobenius map $\thh(R)^\wedge_p\to \thh(R)^{tC_p}$ is the connective cover map. In addition, since both have $\pi_{-1}=0$, the above map is an isomorphism on $\tau_{\geq -1}$ covers. The functor $(-)^{hG}$ for any Lie group $G$ preserves coconnectivity, that is: it takes maps that induce isomorphisms on homotopy groups above a given dimension to maps that induce isomorphisms on homotopy groups above that dimension.  So we apply this to $\phi_p^{hC_{p^n}}:\left(\thh(R)^\wedge_p\right)^{hC_{p^n}}\to \thh(R)^{tC_{p^{n+1}}}$ to see that it is also an equivalence on $\tau_{\geq -1}$-covers. Furthermore, the Frobenius is an $\mathbb{E}_\infty$-ring map, and so we get the claimed ring structure on $\pi_*\left(\tau_{\geq -1}\thh(R)^{hC_{p^n}}\right)$. In addition, we may use this presentation to take the Frobenius to be the isomorphism sending $u\mapsto\sigma$ on homotopy groups. 

It remains to study what the canonical map does. Consider the following commutative diagram:

\[
\begin{tikzcd}
\tau_{\geq -1}\TC^{-}(R)^\wedge_p \arrow[rr, "can^{h\TT}"] \arrow[d]                          &  & \tau_{\geq -1}\operatorname{TP}(R)^\wedge_p \arrow[d] \\
\tau_{\geq -1}\left(\thh(R)^\wedge_p\right)^{hC_{p^n}} \arrow[rr, "can^{hC_{p^n}}"] &  &\tau_{\geq -1} \thh(R)^{tC_{p^{n}}}                         
\end{tikzcd}\]
Applying $\pi_*$, it gives the commutative diagram
\[
\begin{tikzcd}
A_{\textrm{inf}}[u] \arrow[r, "u\mapsto \xi \sigma"]  \arrow[d]   & A_{\textrm{inf}}[\sigma] \arrow[d, "\Tilde{\theta}_{n}"]\\
W_{n+1}(R)[u] \arrow[r] & W_{n}(R)[\sigma]
\end{tikzcd}
\]
The upper horizontal map is identified by \parencite[Proposition 6.3]{BMS2}, and the right vertical map is identified by \parencite[Remark 6.6]{BMS2}. Since we are using the presentation of $\pi_{\geq -1}\left(\left(\thh(R)^\wedge_p\right)^{hC_{p^n}}\right)$ coming from the equivalence with the Tate construction above, it follows that the left vertical map in the above diagram must be $\Tilde{\theta}_{n+1}\circ \phi$. All the maps are surjective on $\pi_0$, so there is a unique map making the diagram commute. From \parencite[Lemma 3.4]{BMS1}, we see that the map on $\pi_0$ must be the restriction, and then the polynomial generator does as claimed.\qed

Since there are suspensions in Equation~(\ref{eqn: Wedge decomp of THH}), we need to understand what happens on the negative homotopy groups as well.
\begin{prop}\label{prop: negative homotopy groups}
The map\[\tau_{< 0}\left(\thh(R)^\wedge_p\right)^{hC_{p^n}}\xrightarrow{can^{hC_{p^{n-1}}}}\tau_{< 0}\left(\thh(R)^{tC_p}\right)^{hC_{p^{n-1}}}\simeq \tau_{< 0}\thh(R)^{tC_{p^{n}}}\] induces an isomorphism
on homotopy groups.
\end{prop}

\pf By the Tate orbit lemma \parencite[Lemma I.2.1]{Nikolaus_Scholze} and induction, as in \parencite[Lemma II.4.1]{Nikolaus_Scholze}, we see that $\left(\left(\thh(R)^\wedge_p\right)_{hC_{p}}\right)^{hC_{p^{n-1}}}\simeq \left(\thh(R)^\wedge_p\right)_{hC_{p^n}}$. Homotopy orbits for any compact Lie group preserve connectivity, so $\left(\thh(R)^\wedge_p\right)_{hC_{p^n}}$ is connective. On the other hand, $\left(\left(\thh(R)^\wedge_p\right)_{hC_{p}}\right)^{hC_{p^{n-1}}}$ is the homotopy fiber of the map $can^{hC_{p^{n-1}}}$. It follows that this map is an equivalence below degree zero, as claimed.\qed

\begin{cor}\label{cor: Complete homotopy Cp calculation}
The map $\TC^{-}(R)\to \left(\thh(R)^\wedge_p\right)^{hC_{p^n}}$ on homotopy groups is the quotient map\[A_{\textrm{inf}}[u,v]/(uv-\xi)\to A_{\textrm{inf}}[u,v]/(uv-\xi, \Tilde{\xi}_nv)\]
\end{cor}
\pf From the computations in Proposition~\ref{prop: frob and can computation} and Proposition~\ref{prop: negative homotopy groups}, we see that the map $\TC^{-}(R)^\wedge_p\to \left(\thh(R)^\wedge_p\right)^{hC_{p^n}}$ must be surjective on homotopy groups. Further, the element $\Tilde{\xi}_nv\in \pi_{-2}\left(\TC^{-}(R)^\wedge_p\right)$ must be mapped to zero in $\pi_{-2}\left(\left(\thh(R)^\wedge_p\right)^{hC_p^n}\right)$ since it goes to zero under the canonical map, which is an isomorphism in this degree. Consequently, we get a factorization \begin{center}
    \begin{tikzcd}
{A_{\textrm{inf}}[u,v]/(uv-\xi)} \arrow[rr] \arrow[d]                &  & \pi_*\left(\thh(R)^\wedge_p\right)^{hC_{p^n}} \\
{A_{\textrm{inf}}[u,v]/(uv-\xi, \Tilde{\xi}_nv)} \arrow[rru, dashed] &  &                                              
\end{tikzcd}
\end{center}
Composing the dotted map with $\phi^{hC_{p^n}}$ is then an isomorphism in non-negative degrees by the computation in Proposition~\ref{prop: hmtpy groups of T fixed points}. Since $\phi_p^{hC_{p^n}}$ is an isomorphism in these degrees by Proposition~\ref{prop: frob and can computation}, it follows that the dotted map must also be an isomorphism.

Similarly, composing with $can^{hC_{p^n}}$ is an isomorphism on negative homotopy groups. Therefore, the dotted map is also an isomorphism on negative homotopy groups by Proposition~\ref{prop: negative homotopy groups}.
\qed

Note that this description, when combined with the isomorphism $\pi_*\left(\thh(R)^{tC_{p^n}}\right)\cong A_{\textrm{inf}}[\sigma,\sigma^{-1}]/(\Tilde{\xi}_n\sigma^{-1})$, completely determine what the Frobenius and the canonical maps must be in this presentation.

\subsection{\texorpdfstring{Case One: $s$ even}{Case One: s even}}
Using the above, we are now ready to evaluate $\tp(A_d, R)^\wedge_p$ and $\TC^{-}(A_d,R)$. By the decomposition in Equation~(\ref{eqn: product decomp of THH}), we may break up these computations into the cases in Equation~(\ref{eqn: Wedge decomp of THH}). In this subsection, we will compute these invariants for the terms with $s$ even.

\begin{lem}\label{lem: sevenresult}
Let $X(m)=\prod_{\substack{s\mid m\\ s\textrm{ even}}}\ \prod_{\overline{\omega}\in \omega_{s,d}}\left(\Sigma^{m-1}\thh(R)\wedge \left(\TT/C_{\frac{m}{s}}\right)_+\right)$, so that we have $X=\prod_{m\in 2\ZZ_+}X(m)$ is the summand of $s$ even elements in $\thh(A_d, R)$. Then there are isomorphisms \[\pi_*\left(\TC^{-}(X(m)^\wedge_p)\right)\cong 
\begin{cases}
\prod\limits_{\substack{s\mid m\\ s\textrm{ even}}}\prod\limits_{\overline{\omega}\in \omega_{s,d}} W_{v_p\left(\frac{m}{s}\right)+1}(R)  & \textrm{if }*> m-1\textrm{ is even}\\
\prod\limits_{\substack{s\mid m\\ s\textrm{ even}}}\prod\limits_{\overline{\omega}\in \omega_{s,d}} W_{v_p\left(\frac{m}{s}\right)}(R)  & \textrm{if }*< m-1\textrm{ is even}\\
0   &   \textrm{otherwise}
\end{cases}\]
and 
\[\pi_{*}\left(\tp(X(m))^\wedge_p\right)\cong \begin{cases}
\prod\limits_{\substack{s\mid m\\ s\textrm{ even}}}\prod\limits_{\overline{\omega}\in \omega_{s,d}}W_{v_p\left(\frac{m}{s}\right)}(R)  &\textrm{if }*\textrm{ is even}\\
0   &\textrm{otherwise}
\end{cases}\]
where $W_0(R)$ is understood to be zero, and $v_p$ is the $p$-adic valuation.
\end{lem}
\pf Since both the homotopy fixed points and the Tate construction will commute with the products in our decomposition, we may compute each term wise. We will compute this for homotopy fixed points, the Tate computation will follow similarly. Thus we need only compute \[\left(\left(\Sigma^{m-1}\thh(R)\wedge \left(\TT/C_{\frac{m}{s}}\right)_+\right)^\wedge_p\right)^{h\TT}\]

For any spectrum $X\in \Sp^{B\TT}$, the map $X\wedge\left(\TT/C_{\frac{m}{s}}\right)_+\xrightarrow{(p\textrm{-completion})\wedge id}X^\wedge_p\wedge \left(\TT/C_{\frac{m}{s}}\right)_+$ mod $p$ becomes up to equivalence $X/p\wedge \left(\TT/C_{\frac{m}{s}}\right)_+\xrightarrow{id} X/p\wedge \left(\TT/C_{\frac{m}{s}}\right)_+$. Consequently, $\left(X\wedge\left(\TT/C_{\frac{m}{s}}\right)_+\right)^\wedge_p\simeq \left(X^\wedge_p\wedge\left(\TT/C_{\frac{m}{s}}\right)_+\right)^\wedge_p$. This equivalence is equivariant, so we have that the above computation reduces to computing \[\left(\left(\Sigma^{m-1}\thh(R)^\wedge_p\wedge \left(\TT/C_{\frac{m}{s}}\right)_+\right)^\wedge_p\right)^{h\TT}\simeq \Sigma^{m-1}\left(\left(\thh(R)^\wedge_p\wedge \left(\TT/C_{\frac{m}{s}}\right)_+\right)^{h\TT}\right)^\wedge_p.\] (Here, to get the $p$-completion inside the Tate construction, we would use \parencite[Lemma I.2.9]{Nikolaus_Scholze}, along with the p-adic equivalence $BC_{p^\infty}\to B\TT$.). By Proposition~\ref{prop: Induced/coinduced action and fixed points}, we in turn get the equivalence 
\begin{align*}
\Sigma^{m-1}\left(\thh(R)^\wedge_p\wedge \left(\TT/C_{\frac{m}{s}}\right)_+\right)^{h\TT}    &\simeq \Sigma^{m-1}\left(S^1\wedge \thh(R)^\wedge_p\right)^{hC_{\frac{m}{s}}}\\
&\simeq_{p\textrm{-equiv}} \Sigma^m \left(\thh(R)^\wedge_p\right)^{hC_{p^{v_p\left(\frac{m}{s}\right)}}}.
\end{align*}
The last equivalence comes from the fact that the map $C_r\to C_{p^{v_p(r)}}$ is a $p$-adic equivalence. The homotopy groups of the third term are completely described by Corollary~\ref{cor: Complete homotopy Cp calculation}. In particular, the third term is already $p$-complete.
\qed

Note that both $\TC^{-}(X)$ and $\tp(X)$ commute with the product. Thus the computation above also lead to a computation of $\TC^-(X)^\wedge_p$ and $\tp(X)^\wedge_p$.

\begin{rem}
It is easy to check that the maps in Lemma~\ref{lem: sevenresult} are isomorphisms of graded $A_{\textrm{inf}}$-modules. Since $\TC^{-}(A_d)$ and $\tp(A_d)$ are $\mathbb{E}_\infty$-ring spectra and $X$ is a summand of $\thh(A_d)$, the graded $A_{\textrm{inf}}$-modules above should have a product. It would be interesting to find out what that product is. This seems tractable, but it is not needed for the results of this paper.
\end{rem}

\subsection{\texorpdfstring{Case Two: $m$ odd}{Case Two: m odd}} In this section, we compute the topological negative cyclic and periodic homology in the the second case in Equation~(\ref{eqn: Wedge decomp of THH}).

\begin{lem}\label{soddresult}
Let $X(m)=\prod_{s\mid m}\prod_{\overline{\omega}\in \omega_{s,d}} \left(\Sigma^{m-1}\thh(R)\wedge \left(\TT/C_{\frac{m}{s}}\right)_+\right)$, so that if we let $X=\prod_{m\in \ZZ_+\setminus 2\ZZ_+} X(m)$, then $X$ is the summand where $s$ and $m$ are odd in Equation~(\ref{eqn: Wedge decomp of THH}). Then there are isomorphisms
\[
\pi_*\left(\TC^{-}(X(m))^\wedge_p\right)\cong 
\begin{cases}
\prod\limits_{s\mid m}\prod\limits_{\overline{\omega}\in \omega_{s,d}}W_{v_p\left(\frac{m}{s}\right)+1}(R)    &   \textrm{if }*> m-1\textrm{ is odd}\\
\prod\limits_{s\mid m}\prod\limits_{\overline{\omega}\in \omega_{s,d}}W_{v_p\left(\frac{m}{s}\right)}(R)    &   \textrm{if }*< m-1\textrm{ is odd}\\
0   &   \textrm{otherwise}
\end{cases}
\]
and
\[
\pi_*\left(\tp(X(m))^\wedge_p\right)\cong\begin{cases}
\prod\limits_{s\mid m}\prod\limits_{\overline{\omega}\in \omega_{s,d}}W_{v_p\left(\frac{m}{s}\right)}(R)  &   \textrm{if }*\textrm{ is odd}\\
0   &   \textrm{otherwise}
\end{cases}
\]
where $W_0(R)$ is understood to be zero, and $v_p$ is the $p$-adic valuation. 
\end{lem}\label{lem: case two hmtpy and tate}
\pf There was nothing essential about the parity of $m$ and $s$ used in the proof of Lemma~\ref{lem: sevenresult}. Hence the same proof will go through in this case, the only difference being the parity of $m$ is odd, so the homotopy groups will be concentrated in odd degree instead of even degree.\qed

\subsection{\texorpdfstring{Case Three: $m\neq s\mod 2$}{Case Three: m odd and s even}}

We will now deal with the third case for $R$ a perfectoid ring with respect to the prime $p=2$.

\begin{lem}\label{lem: m neq s mod 2}
Let $X(m)=\prod_{\substack{s\mid m \\ s\neq m\mod 2}}\prod_{\overline{\omega}\in w_{s,d}}\left(\Sigma^{m-2}\thh(R)\wedge \RP^2\left(\frac{m}{s}\right)\right)$, so that $X:=\prod_{m\in 2\ZZ}X(m)$ is the summand of Equation~\ref{eqn: product decomp of THH} of elements with $m\neq s\mod 2$. Then there are equivalences \[\pi_*\left(\tp(X(m))^\wedge_2\right)\simeq\pi_*\left(\TC^{-}(X(m))^\wedge_2\right)\simeq \begin{cases}
\prod\limits_{\substack{s\mid m \\ s\neq m\mod 2}}\prod\limits_{\overline{\omega}\in \omega_{s,d}}R & \textrm{If }*\textrm{ is odd}\\
0 & \textrm{otherwise}
\end{cases}\]

\end{lem}

In order to prove this lemma, we must make a few reductions. By the same argument as in Lemma~\ref{lem: sevenresult}, we may consider the homotopy and Tate fixed point construction one factor at a time.  Recall that we have a cofiber sequence \[\left(\TT/C_{\frac{m}{2s}}\right)_+\to \left(\TT/C_{\frac{m}{s}}\right)_+\to \RP^2\left(\frac{m}{s}\right)\] from which we get the cofiber sequences \[\left(\Sigma^{m-1}\thh(R)^{hC_{\frac{m}{2s}}}\right)^\wedge_2\to \left(\Sigma^{m-1}\thh(R)^{hC_{\frac{m}{s}}}\right)^\wedge_2\to \left(\left(\Sigma^{m-2}\thh(R)\wedge \RP^2\left(\frac{m}{s}\right)\right)^{h\TT}\right)^\wedge_2\]
and 
\[\left(\Sigma^{m-1}\thh(R)^{tC_{\frac{m}{2s}}}\right)^\wedge_2\to \left(\Sigma^{m-1}\thh(R)^{tC_{\frac{m}{s}}}\right)^\wedge_2\to \left(\left(\Sigma^{m-2}\thh(R)\wedge \RP^2\left(\frac{m}{s}\right)\right)^{t\TT}\right)^\wedge_2.\] 
Thus in order to compute the homotopy and Tate fixed point construction on $X(m)$ we will first need to identify the maps induced by $\TT/C_{\frac{m}{2s}}\to \TT/C_{\frac{m}{s}}$ on homotopy groups. In fact, these maps are the Verscherbung maps, which can be seen from the following lemma:

\begin{lem}\label{lem: Versch map}
Let $R$ be a perfectoid ring with respect to a prime $p$. Then the transition maps \[\left(\thh(R)^\wedge_p\right)^{hC_{p^{n-1}}}\to \left(\thh(R)^\wedge_p\right)^{hC_{p^{n}}}\]
and 
\[\left(\thh(R)^\wedge_p\right)^{tC_{p^{n-1}}}\to \left(\thh(R)^\wedge_p\right)^{tC_{p^{n}}}\] on homotopy groups are the Verschiebung maps $V:W_{n-1}(R)\to W_{n}(R)$, or $V:W_{n}(R)\to W_{n+1}(R)$ in the case of the homotopy fixed points in even positive degrees. 
\end{lem}

\pf Note that since these maps are induced by the map of $\thh(R)^\wedge_p$-modules $\thh(R)^\wedge_p\wedge \left(\TT/C_{\frac{m}{2s}}\right)_+\to \thh(R)^\wedge_p\wedge \left(\TT/C_{\frac{m}{s}}\right)_+$, it follows that the maps in question on homotopy groups are $A_{\textrm{inf}}[u,v]/(uv-\xi)$-module and $A_{\textrm{inf}}[\sigma^{\pm 1}]$-module maps, respectively. In particular, since both are principal modules over their respective rings by Corollary~\ref{cor: Complete homotopy Cp calculation}, we need only determine what the map does on $\pi_0$. 

In order to identify this map, we must use particular features of the identification of $\pi_0\left(\thh(R)^{tC_{p^n}}\right)$ as $W_{n}(R)$. We record here the identification of Bhatt-Morrow-Scholze in \cite[Remark 6.6]{BMS2}. We first recall that we have the isotropy separation pullback square 
\[\begin{tikzcd}
\left(\thh(R)^\wedge_p\right)^{C_p} \arrow[r] \arrow[d] & \left(\thh(R)^\wedge_p\right)^{\Phi C_p} \arrow[d] \\
\left(\thh(R)^\wedge_p\right)^{hC_p} \arrow[r]          & \left(\thh(R)^\wedge_p\right)^{tC_p}              
\end{tikzcd}\] where $(-)^{C_p}$ and $(-)^{\Phi C_p}$ are the genuine and geometric fixed points functors, respectively.

In the framework of genuine cyclotomic spectra, which $\thh$ fits into, we then also have an equivalence $\thh(R)^\wedge_p\simeq \rho_p^*\left(\thh(R)^\wedge_p\right)^{\Phi C_p}$, $\rho_p:\TT\to\TT/C_p$ a homeomorphism, such that the composition \[\thh(R)^\wedge_p \simeq \rho_p^*\left(\thh(R)^\wedge_p\right)^{\Phi C_p}\to \varrho_p^*\left(\thh(R)^\wedge_p\right)^{tC_p}\] is the Frobenius. In our case by Proposition~\ref{prop: frob and can computation} the composition is then an isomorphism above degree $-1$, and so by general properties of pullback squares the map $\left(\thh(R)^\wedge_p\right)^{C_p}\to \left(\thh(R)^\wedge_p\right)^{hC_p}$ is an isomorphism upon taking connective covers. By the main result of \cite{Tsalidis}, we then have that the maps $\left(\thh(R)^\wedge_p\right)^{C_{p^n}}\to \left(\thh(R)^\wedge_p\right)^{hC_{p^n}}$ are isomorphisms upon taking connective covers for $n\geq 1$.

Thus we have isomorphisms \[W_{n+1}(R)\simeq \operatorname{TR}^n_0(R)^\wedge_p\simeq \pi_0\left(\thh(R)^\wedge_p\right)^{hC_{p^n}}\simeq \pi_0\left(\thh(R)^\wedge_p\right)^{tC_{p^{n+1}}}\]
the first equivalence coming from \cite[Theorem F]{Hesselholt_Madsen}. From this discussion, it also follows that the transition maps on $\pi_0$ are then also the maps $\operatorname{TR}^{n-1}_0(R)\to \operatorname{TR}^{n}_0(R)$, which we may identify as the Verschiebung map by \cite[Theorem 3.3]{Hesselholt_Madsen}.\qed

We are now ready to prove Lemma~\ref{lem: m neq s mod 2}.

\textbf{Proof of Lemma~\ref{lem: m neq s mod 2}. \ } We will prove this for the Tate fixed points, the homotopy fixed points being similar. Write $X(m,s,\overline{\omega})$ for the factor of $X(m)$ corresponding to $s$ and $\overline{\omega}$\footnote{This is essentially $\thh(R)\wedge B^{cy}(\overline{\omega})$, but we use this notation instead to emphasize that the representation sphere has been trivialized.}. Then the cofiber sequence \[\left(\Sigma^{m-1}\thh(R)^{tC_{\frac{m}{2s}}}\right)^\wedge_2\to \left(\Sigma^{m-1}\thh(R)^{tC_{\frac{m}{s}}}\right)^\wedge_2\to \left(\left(\Sigma^{m-2}\thh(R)\wedge \RP^2\left(\frac{m}{s}\right)\right)^{t\TT}\right)^\wedge_2\] gives, since the homotopy groups of both $\left(\Sigma^{m-1}\thh(R)^{tC_{\frac{m}{2s}}}\right)^\wedge_2$ and $\left(\Sigma^{m-1}\thh(R)^{tC_{\frac{m}{s}}}\right)^\wedge_2$ are concentrated in odd degrees, that \[\pi_{*}\left(\left(X(m,s,\overline{\omega})^\wedge_2                                                        \right)^{t\TT}\right)\cong \begin{cases}
\ker\left(\pi_{*-2}\left(\thh(R)^{tC_{\frac{m}{2s}}}\right)^\wedge_2\xrightarrow{V}\pi_{*-2}\left(\thh(R)^{tC_{\frac{m}{s}}}\right)^\wedge_2\right)   &   \textrm{If }*\textrm{ even}\\
\coker\left(\pi_{*-1}\left(\thh(R)^{tC_{\frac{m}{2s}}}\right)^\wedge_2\xrightarrow{V}\pi_{*-1}\left(\thh(R)^{tC_{\frac{m}{s}}}\right)^\wedge_2\right) &\textrm{otherwise}
\end{cases}\] and then the above lemma gives the result.\qed

\section{\texorpdfstring{$\TC$ calculation}{TC calculation}}\label{sec: tc}
We can now calculate $\pi_*(\TC(A_d)^\wedge_p)$. We start with the case $p\neq 2$, and will handle the case $p=2$ in Subsection~\ref{sec: p=2 final computation}.  To do this, we will use the technique introduced by Nikolaus and Scholze.
\begin{thm}[\cite{Nikolaus_Scholze}, Theorem II.4.11] Let $X$ be a genuine cyclotomic spectrum such that the underlying spectrum is bounded below. Then 
\[\TC(X)\simeq\hofib\left(\TC^{-}(X)\xrightarrow{(\phi_p-can)_{p\in \mathbb{P}}}\prod_{p\in \mathbb{P}}\tp(X)^\wedge_p\right). \]  
\end{thm}

Taking the $p$-completion of the above diagram, or using the $p$-adic equivalence $BC_{p^\infty}\to B\TT$ and \parencite[Theorem II.4.10]{Nikolaus_Scholze}, we get the following
\begin{align*}
\TC(X)^\wedge_p &\simeq \hofib\left(\TC^{-}(X)^\wedge_p\xrightarrow[]{\phi_p-can}\tp(X)^\wedge_p\right)\\
                &\simeq \hofib\left(\TC^{-}(X^\wedge_p)\xrightarrow[]{\phi_p-can} \tp(X^\wedge_p)^\wedge_p\right)
\end{align*}
For
$X=\thh(A_d)$, the splitting of $\thh(A_d)$ in Equation (\ref{eqn: Wedge decomp of THH}) can be written as a product as explained below Equation (\ref{eqn: product decomp of THH})  and since we are $p$-completing and $p\neq 2$, $(\RP ^2)^\wedge_p\simeq *$ and so the factors with $s$ odd and $m$ even collapse.  This gives a $\TT$-invariant decomposition
\begin{equation}
\thh(A_d)^\wedge_p \simeq \prod_{m\in \NN}\quad \prod_{\substack{s\mid m\\ s\equiv m\textrm{ mod}\ 2}}
\quad \prod_{\overline{\omega}\in \omega_{s,d}} \left(\Sigma^{m-1}\thh(R) \wedge \left(\TT/C_{\frac{m}{s}}\right)_+\right)^\wedge_p
\label{eqn: Prod decomp of THH}
\end{equation}
with actions as explained in the previous sections.

We know what $\TC^-$ and $\tp$ look like on each of the factors: Lemma \ref{lem: sevenresult} deals with the factors where $s$ and $m$ are both even and Lemma \ref{soddresult} deals with the factors where they are both odd.

We can determine what $\phi_p$ and the canonical map $can$ do in these terms.  Note that by Lemma \ref{lem:smash}, the Frobenius on $\thh(R)\wedge B^{cy}(\Pi) \simeq \thh(R[\Pi])$ is induced by the smash product of the Frobenius maps on each part.  The decomposition of Equation~(\ref{Bsplit}) is one of $\TT$-spaces, and underlies the decompositions of $\TC^-$ and $\tp$ that we are using.  The canonical map goes from the homotopy fixed points to the homotopy Tate construction of a particular $\TT$-space, so it sends the summand corresponding to a particular $\overline{\omega}$ (which determines its $m$ and $s$) to itself.  However, the Frobenius map on $B^{cy}(\Pi)$ sends the summand corresponding to  $\overline{\omega}$ (of length $m$ and cycle length $s$) to the concatenation $\overline{\omega ^{\star p}}$ of $\overline{\omega}$ $p$ times (of length $pm$ and with the same cycle length
$s$).

Before considering the two cases, we identify the Frobenius in terms of the Frobenius on $\thh(R)^\wedge_p$. Indeed, we have the following commutative diagram:
\[
\begin{tikzcd}
\thh(R)^\wedge_p\wedge B^{cy}(\overline{\omega}) \arrow[rr, "\phi_p^{\thh(R)}\wedge id"] \arrow[rrdd, bend right, "\phi_p"'] &  & \thh(R)^{tC_p}\wedge B^{cy}(\overline{\omega}) \arrow[d, "id\wedge \phi_p^{B}"] \\
                                                                                                                         &  & \thh(R)^{tC_p}\wedge B^{cy}(\overline{\omega^{\ast p}})^{tC_p} \arrow[d, "l"]           \\
                                                                                                                         &  & \left(\thh(R)^\wedge_p\wedge B^{cy}(\overline{\omega^{\ast p}})\right)^{tC_p}          
\end{tikzcd}
\]
where $l$ is the lax-monoidal map. Then, by a similar argument as in \parencite[Lemma 2]{Hesselholt_Nikolaus}, the vertical maps in the above diagram are equivalences. Thus, up to a canonical isomorphism, we may take $\phi_p=\phi_p^{\thh(R)}\wedge id$.

\subsection{\texorpdfstring{Case One: $p> 2$}{Case One: p>2}}
We consider the case of $s$ even first. By Lemmas \ref{lem: sevenresult} and \ref{soddresult}, the homotopy groups of both target and source are concentrated in even degrees, so we have
\[\pi_{2r} (\TC(X)^\wedge_p) \cong\ker \left(\pi_{2r}\left(\TC^{-}(X)^\wedge_p\right)\xrightarrow{(\phi_p-can)_*}
\pi_{2r}\left(\tp(X)^\wedge_p\right)\right)\]
and
\[\pi_{2r-1} (\TC(X)^\wedge_p) \cong\coker \left(\pi_{2r}\left(\TC^{-}(X)^\wedge_p\right)\xrightarrow{(\phi_p-can)_*}
\pi_{2r}\left(\tp(X)^\wedge_p\right)\right).\]
Where $X=\prod_{m\in 2\ZZ_+}\prod_{\substack{s\mid m\\ s\textrm{ even}}}\prod_{\overline{\omega}\in \omega_{s,d}}\left(\Sigma^{m-1}\thh(R)\wedge \left(\TT/C_{\frac{m}{s}}\right)_+\right)$. Since the Frobenius fixes $s$, and only changes $m$, we will rewrite \[X = \prod_{m'\in 2J_p}\prod_{\nu = 0}^\infty \prod_{\substack{s\mid m'p^\nu\\ s \textrm{ even}}}\prod_{\overline{\omega}\in \omega_{s,d}}\Sigma^{m-1}\thh(R)\wedge \left(\mathbb{T}/C_{\frac{m}{s}}\right)_+\]
where $J_p$ is the set of positive integers coprime to $p$ (as in \cite{Speirs_coordinate_axes}). Both $\phi_p$ and the canonical map $can$ will respect the outside product, so it is enough to consider $X(m')$ for fixed $m'\in 2J_p$, $X(m')$ being the $m'$ term in $X$. For the rest of this subsection, we will have $m'$, $p$, and $r$ fixed.

We see from Lemma~\ref{lem: sevenresult} that  there are two distinct cases to consider. We will split them up similar to the argument in \cite{Hesselholt_Nikolaus}. Recall that the function $t_{ev}=t_{ev}(p, r, m')$ is the unique positive integer, if it exists, such that $m'p^{t_{ev}-1}\leq 2r<m'p^{t_{ev}}$. If no such positive integer exists, $t_{ev}(p,r,m')=0$.  Then we have the following commutative diagram, where the columns are the exact sequences:

\[
\begin{tikzcd}
0 \arrow[d]                                                                                                                                                                                    &  & 0 \arrow[d]                                                                                                                                     \\
{\prod\limits_{\nu = t_{ev}}^\infty\ \prod\limits_{\substack{s\mid m'p^\nu\\ s\textrm{ is even}}}\ \prod\limits_{\omega_{s,d}}W_{\nu - v_p(s)}(R)} \arrow[d] \arrow[rr, "\phi_p-can"]               &  & {\prod\limits_{\nu = t_{ev}}^\infty\ \prod\limits_{\substack{s\mid m'p^\nu\\ s\textrm{ even}}}\ \prod\limits_{\omega_{s,d}}W_{\nu - v_p(s)}(R)} \arrow[d]  \\
\pi_{2r}\left(\TC^{-}(X(m'))^\wedge_p\right) \arrow[d] \arrow[rr, "\phi_p-can"]                                                                                                                                &  & \pi_{2r}\left(\tp(X(m'))^\wedge_p\right) \arrow[d]                                                                                                                       \\
{\prod\limits_{\nu = 0}^{t_{ev}-1}\ \prod\limits_{\substack{s\mid m'p^\nu\\ s\textrm{ even}}}\ \prod\limits_{\omega_{s,d}}W_{\nu - v_p(s)+1}(R)} \arrow[d] \arrow[rr, "\overline{\phi_p-can}"] &  & {\prod\limits_{\nu = 0}^{ t_{ev}-1}\ \prod\limits_{\substack{s\mid m'p^\nu\\ s\textrm{ even}}}\ \prod\limits_{\omega_{s,d}}W_{\nu - v_p(s)}(R)} \arrow[d] \\
0                                                                                                                                                                                              &  & 0                                                                                                                                              
\end{tikzcd}
\]

Recall that by Proposition~\ref{prop: frob and can computation} and Proposition~\ref{prop: negative homotopy groups}, we know what the Frobenius and the canonical map look like on each component, at least enough to do computations. On the top line, the canonical map is an isomorphism. In the bottom components, the Frobenius can be taken to be the identity, and the canonical map is the reduction in length and multiplication by $\Tilde{\theta}_r(\xi)$ on the polynomial generator.

\begin{lem}\label{lem: reduction to good case}
The top horizontal map in the above diagram is an isomorphism.
\end{lem}
\pf
In this range, the canonical map is an isomorphism. Let \[x\in \prod_{\nu = t_{ev}}^\infty\prod_{\substack{s\mid m'p^\nu\\ s\textrm{ even}}}\prod_{\omega_{s,d}}W_{\nu - v_p(s)}(R),\] and let $p_{t_{ev}}, p_{t_{ev}+2}, \ldots$ be the projection maps of the outermost product. Define the degree of $x$ to be the minimum $n$ such that $p_n(x)\neq 0$, which exists if and only if $x\neq 0$. Then $can(x)$ also has degree $n$ since $can$ is an isomorphism on the product terms. On the other hand, $\phi_p$ must increase the degree by at least one, or send $x$ to zero. In either case, $p_n(\phi_p(x))=0$, so $p_n(can(x)-\phi_p(x))\neq 0$, and $can-\phi_p$ is injective. 

To see that $can-\phi_p$ is surjective, we write $x_n := p_n(x)$. We will define a pre-image $y$ inductively as $y_{t_{ev}}=can^{-1}(x_{t_{ev}})$, and $y_k = can^{-1}(x_k - \phi_p(y_{k-1}))$. It is then straightforward to check that $(can-\phi_p)(y)=x$, as desired.\qed

Hence we get by the snake lemma applied to the diagram above that $\pi_{2r}\left(\TC(X(m'))^\wedge_p\right)\cong \ker(\overline{can-\phi_p})$, and $\pi_{2r}\left(\TC_{2r-1}(X(m'))^\wedge_p\right) \cong \coker(\overline{can-\phi_p})$. Once again, to compute $\overline{can-\phi_p}$, there are two cases to consider: when $\nu = t_{ev}-1$ and when it does not.

For convenience, we write \[\TC^{-}_{2r}(+):=\prod\limits_{\nu = 0}^{t_{ev}-1}\ \prod\limits_{\substack{s\mid m'p^\nu\\ s \textrm{ even}}}\ \prod\limits_{\overline{\omega}\in \omega_{s,d}}W_{\nu - v_p(s)+1}(R)\] and \[\tp_{2r}(+):=\prod\limits_{\nu = 0}^{t_{ev}-1}\ \prod\limits_{\substack{s\mid m'p^\nu\\ s\textrm{ even}}}\ \prod\limits_{\overline{\omega}\in \omega_{s, d}}W_{\nu - v_p(s)}(R).\] We then have the following result.

\begin{lem}\label{lem: almost iso}
The composite map \[\prod_{\nu = 0}^{t_{ev}-2}\prod_{\substack{s\mid m'p^\nu \\ s\textrm{ even}}}\ \prod_{\omega_{s,d}}W_{\nu - v_p(s)+1}(R)\hookrightarrow \TC^{-}_{2r}(+)\xrightarrow[]{\overline{can-\phi_p}}\tp_{2r}(+)\] is an isomorphism.
\end{lem}
\pf In this range, $\overline{can-\phi_p}=can-\phi_p$, the difference only comes from the $\nu = t_{ev}-1$ factor in $\TC^-_{2r}(+)$, where the canonical map is unaffected, but the Frobenius map is the zero map. In addition, $\phi_p$ is an isomorphism on each component of the source. We will now transfer the proof of Lemma~\ref{lem: reduction to good case} to this setting. Let $x\in \prod_{\nu=0}^{t_{ev}-2}\prod_{\substack{s\mid m'p^{\nu}\\ s \textrm{ even}}}\prod_{\omega_{s, d}}W_{\nu - v_p(s)+1}(R)$ be nonzero. Define the degree of $x$ as the \textit{largest} $n$ such that $p_n(x)\neq 0$, $p_i$ the projections of the outermost product. Then $\phi_p(x)$ has degree $n+1$, and $can(x)$ has degree at most $n$. Consequentially, the degree of $(can-\phi_p)(x)$ is $n+1$, and in particular $(can-\phi_p)(x)\neq 0$.

To see that the map is also surjective, we first show that $\phi_p$ is surjective, i.e. it hits all the non-zero factors in the product. Since $\omega_{s,d}$ depends only on $s$ and $d$, and $\phi_p$ preserves both, it follows that $\phi_p$ will hit all factors in this product if it hits one. For the product over $s$, the Frobenius will not change $s$, so we are getting the factors with $s$ such that $s\mid m'p^{\nu}$ in degree $\nu+1$. In other words, the Frobenius is only missing the factors with $v_p(s)=v_p(m)$. These terms on the right, however, are exactly the $W_0$-terms, and so are all zero. Finally, in the outermost product we only miss the factor $\nu=0$ which is also only $W_0$-terms.

Thus, $\phi_p$ is surjective in this range. To see that $\phi_p-can$ is surjective, take an element $y\in \prod_{\nu=0}^{\nu=t_{ev}-1}\prod_{\substack{s\mid m'p^\nu\\ s \textrm{ even}}}\prod_{\overline{\omega}\in \omega_{s,d}}W_{\nu - v_p(s)}(R)$. Inductively define a pre-image $z$ by $z_{t_{ev}-2}=\phi_p^{-1}(y_{t_{ev}-1})$, and $z_{t_{ev}-k} = \phi_p^{-1}(y_{t_{ev}-k}-can(z_{t_{ev}-k+1}))$.\qed

Thus the cokernel of the map in question is zero. Furthermore, this also shows that the kernel is, in addition to $0$, comprised of elements $X$ such that $x_{t_{ev}-1}\neq 0$. For any such $x$, for it to be in the kernel, $x_{t_{ev}-2}=-\phi_p^{-1}(can(x_{t_{ev}}))$. Inductively, $x_i$ is determined by $x_{t_{ev}-1}$. Conversely, every element in $\prod_{\substack{s\mid m'p^{t_{ev}-1}\\s \textrm{ even}}}\prod_{\omega_{s,d}}W_{t_{ev}-v_p(s)}(R)$ determines an element in the kernel, so $p_{t_{ev}-1}|_{\ker(\overline{can-\phi_p})}$ is an isomorphism.
\begin{cor}
Let $X$ be as above, p an odd prime. Then there are isomorphisms \[\pi_{2r-1}\left(\TC(X)\right)^\wedge_p\cong 0\] and \[\pi_{2r}\left(\TC(X)^\wedge_p\right)\cong \prod_{m'\in 2J_p}\prod_{\substack{s\mid m'p^{t_{ev}-1}\\ s\textrm{ even}}}\prod_{\omega_{s,d}}W_{t_{ev}-v_p(s)}(R).\]
\end{cor}\label{cor: Case One computation}
\pf By the decomposition explained above, it is enough to show that this holds when $m'$ is fixed. This follows from Lemma~\ref{lem: almost iso} and the discussion after it.\qed

The case when $m$ is odd is essentially identical in proof to Case One. For this reason, we will only record the result

\begin{lem}
Let $X=\prod_{m\in \ZZ_+\setminus 2\ZZ_+}\prod_{s\mid m}\prod_{\omega_{s,d}} \left(\Sigma^{m-1}\thh(R)\wedge \left(\TT/C_{\frac{m}{s}}\right)_+\right)$, and let $p$ be an odd prime. Then there are isomorphisms $\pi_{2r}\left(\TC(X)^\wedge_p\right)=0$ and \[\pi_{2r+1}\left(\TC(X)^\wedge_p\right)\cong \prod_{m'\in J_p\setminus 2J_p}\prod_{s\mid m'p^{t_{od}-1}}\prod_{\omega_{s,d}}W_{t_{od}-v_p(s)}(R)\]
\end{lem}

\subsection{\texorpdfstring{Case Two: $p=2$}{Case Two: p=2}}\label{sec: p=2 final computation}

For $p=2$, we still get a decomposition, but for $p=2$ the decomposition is in terms of $s$: we have the case when $s$ is even and when $s$ is odd. For $s$ even, there is no distinction between when $p=2$ and when $p>2$, we may use the same arguments and get the same answer. 

Hence we need only consider the case of $s$ odd. To this end, let \[X(m'):=\prod_{s\mid m'}\ \prod_{\omega_{s,d}} \Sigma^{m'-1}\thh(R)\wedge \left(\TT/C_{\frac{m'}{s}}\right)_+\times \prod_{s\mid m'}\ \prod_{\nu = 1}^{\infty}\ \prod_{ \omega_{s,d}}\Sigma^{m'2^\nu -2}\thh(R)\wedge \RP^2\left(\frac{m'2^\nu}{s}\right)\] so that $\prod_{m'\in \ZZ\setminus 2\ZZ} X(m')$ is exactly the factors of Equation~\ref{eqn: Prod decomp of THH} with $s$ odd. By Lemma~\ref{lem: m neq s mod 2}, we then have isomorphisms \[\pi_{*}\left(\TC^{-}(X(m'))^\wedge_2\right)\simeq \begin{cases}
\prod\limits_{s\mid m'}\ \prod\limits_{\nu = 0}^{\infty}\ \prod\limits_{\omega_{s,d}} R  &   \textrm{If }*\geq m'-1\textrm{ is odd}\\
\prod\limits_{s\mid m'}\ \prod\limits_{\nu = 1}^{\infty}\ \prod\limits_{\omega_{s,d}} R  &   \textrm{If }*< m'-1\textrm{ is odd}\\
0   &   \textrm{otherwise}
\end{cases}\]
and
\[\pi_{*}\left(\tp(X(m'))^\wedge_2\right)\simeq \begin{cases}
\prod\limits_{s\mid m'}\ \prod\limits_{\nu = 1}^{\infty}\ \prod\limits_{\omega_{s,d}} R  &   \textrm{If }*\textrm{ is odd}\\
0   &   \textrm{otherwise}
\end{cases}\]
and so it remains to identify what the Frobenius and canonical maps are. 

In order to identify what these maps are, note that the cofiber sequence defining $\RP^2(i)$ implies that the canonical map is the map induced by the canonical map on $\thh(R)^\wedge_p\wedge \left(\TT/C_{\frac{m}{s}}\right)_+$. The Frobenius is also induced by this cofiber sequence, using the same trick as in the $p\neq 2$ case in order to use the Frobenius on $\thh(R)$ as the Frobenius on the whole spectrum. In particular, both the Frobenius and the canonical maps are still isomorphisms in the ranges they are above. We may run a similar analysis as in Corollary~\ref{cor: Case One computation} to derive the following.

\begin{cor}
Let $X(m')$ be as above. Then there are isomorphisms 
\[
\pi_{*}\left(\TC(X(m'))^\wedge_2\right)\cong \begin{cases}
\prod\limits_{\nu=0}^{t_{ev}(2, r, m')-1}\prod\limits_{s\mid m'}\prod\limits_{\omega_{s,d}}R &   \textrm{if }*=2r+1\\
0   &   \textrm{if }*=2r
\end{cases}.
\]
\end{cor}

\begin{rem}
We may use in the above Corollary either $t_{ev}$ or $t_{od}$, the difference between the two functions only appears for odd primes. 
\end{rem}

Theorem~\ref{thm: TC main thm}, and therefore Theorem~\ref{thm: main result}, now follows.

\section{Examples}\label{sec: examples}
We conclude this paper with some calculations which can be derived from our results. The first example is a direct corollary, but one which tends to be the main class of computations done for $K$-theory. The second is not an example of a perfectoid ring, but is sufficiently close that we may still apply our results. Finally, the last example recovers recent computations done by Speirs. 

\subsection{\texorpdfstring{The case of $R$ a perfect $\FF_p$-algebra}{The case of R a perfect Fp algebra}}\label{sec: perfect algebra}

As stated above, this is a direct application of the main theorem, since all perfect $\FF_p$-algebras are perfectoid.

\begin{cor}
Let $p$ be an odd prime, and $k$ a perfect $\mathbb{F}_p$-algebra, i.e., an $\FF_p$-algebra where the Frobenius is an isomorphism. Then there is an isomorphism \[\pi_*\left(K(A_d, \mathfrak{m})\right)\cong \begin{cases}
\prod\limits_{m'\in 2J_p}\prod\limits_{\substack{s\mid m'p^{t_{ev}-1}\\  s\textrm{ even}}}\prod\limits_{\omega_{s,d}}W_{t_{ev}(p,r,m')-v_p(s)}(k)   & \textrm{if }*=2r\\
\prod\limits_{m'\in J_p\setminus 2J_p}\prod\limits_{s\mid m'p^{t_{od}-1}}\prod\limits_{\omega_{s,d}}W_{t_{od}(p, r, m')-v_p(s)}(k)   &   \textrm{if }*=2r+1
\end{cases}.\]
\end{cor}
This recovers and slightly extends the calculation Lindenstrauss and McCarthy \cite{LM_Extensions_by_direct_sums} did for perfect fields of positive characteristic. Note that we are not $p$-completing relative $K$-theory in the above corollary. This is because of Proposition~\ref{prop: perfect p complete}.

In particular, we get the following result.
\begin{cor}
Let $p$ be an odd prime, and let $R=\FF_{q}$. Then there is an isomorphism
\[
\pi_*\left(K(A_d)\right)\cong\begin{cases}
\ZZ &   \textrm{if }*=0\\
\prod\limits_{m'\in 2J_p}\prod\limits_{\substack{s\mid m'p^{t_{ev}-1}\\  s\textrm{ even}}}\prod\limits_{\omega_{s,d}}\mathcal{O}_F/p^{t_{ev}(p,r,m')-v_p(s)}   & \textrm{if }*=2r,\textrm{ }*\neq 0\\
\ZZ/(q^{r-1}-1)\times \prod\limits_{m'\in J_p\setminus 2J_p}\prod\limits_{s\mid m'p^{t_{od}-1}}\prod\limits_{\omega_{s,d}}\mathcal{O}_F/p^{t_{od}(p, r, m')-v_p(s)}  &   \textrm{if }*=2r+1
\end{cases}
\]
where $F=\QQ_p(\zeta_{q-1})$ is the unique unramified extension of $\QQ_p$ of degree $\log_p(q)$.

For $p=2$, and $R=\FF_{2^n}$, we have isomorpisms
\[
\pi_{*}\left(K(A_d)\right)\cong \begin{cases}
\ZZ    &    \textrm{if }*=0\\
\prod\limits_{m'\in 2\ZZ}\ \prod\limits_{\substack{s\mid m'2^{t_{ev}-1}\\ s\textrm{ even}}}\ \prod\limits_{\omega_{s,d}}\mathcal{O}_{F}/2^{t_{ev}-v_2(s)} &   \textrm{if }*=2r, r>0\\
\ZZ/(2^{n(r-1)}-1)\times \prod\limits_{m'\in \ZZ\setminus 2\ZZ}\prod\limits_{s\mid m'}\prod\limits_{\nu=0}^{t_{ev}-1}\prod\limits_{\omega_{s,d}}\FF_{2^n}   &   \textrm{if }*=2r+1
\end{cases}
\]
\end{cor}

\subsection{\texorpdfstring{The case of $R=\ZZ_p[\zeta_{p^\infty}]$}{The case of R the cyclotomic integers}}
For this example, note that as stated above $\ZZ_p[\zeta_{p^\infty}]$ is not perfectoid. In particular, every perfectoid ring is $p$-complete and this ring is not. This, however, is the only problem: the ring $\ZZ_p^{cycl}=\ZZ_p[\zeta_{p^\infty}]^\wedge_p$ is perfectoid. In addition, since $\ZZ_p[\zeta_{p^\infty}]$ has bounded $p^{\infty}$-torsion (because it is torsion free), this is also the derived $p$-completion. In particular, \parencite[Lemma 5.2]{Clausen_Mathew_Morrow} applies and \[\thh(\ZZ_p[\zeta_{p^\infty}])^\wedge_p\simeq \thh(\ZZ_p^{cycl})^\wedge_p.\]

Since this map comes from the natural map of rings $\ZZ_p[\zeta_{p^\infty}]\to \ZZ_p^{cycl}$, the map on $\thh$ must be a map of cyclotmic spectra, and we may safely replace $\thh(\ZZ_p[\zeta_{p^\infty}])^\wedge_p$ with $\thh(\ZZ_p^{cycl})^\wedge_p$ in all of the above arguments. Hence we get the following:
\begin{cor}
Let $p$ be an odd prime, and let $R=\ZZ_p[\zeta_{p^\infty}]$. Then there are isomorphisms
\[
\pi_*\left(K(A_d,\mathfrak{m})^\wedge_p\right)\cong \begin{cases}
\prod\limits_{m'\in 2J_p}\prod\limits_{\substack{s\mid m'p^{t_{ev}-1}\\  s\textrm{ even}}}\prod\limits_{\omega_{s,d}}W_{t_{ev}(p,r,m,)-v_p(s)}(\ZZ_p^{cycl})   & \textrm{if }*=2r\\
\prod\limits_{m'\in J_p\setminus 2J_p}\prod\limits_{s\mid m'p^{t_{od}-1}}\prod\limits_{\omega_{s,d}}W_{t_{od}(p, r, m')-v_p(s)}(\ZZ_p^{cycl})   &   \textrm{if }*=2r+1
\end{cases}.
\]

If $p=2$, then there are isomorphisms \[
\pi_*\left(K(A_d,\mathfrak{m})^\wedge_2\right)\cong \begin{cases}
\prod\limits_{m'\in 2\ZZ}\prod\limits_{\substack{s\mid m'2^{t_{ev}-1}\\  s\textrm{ even}}}\prod\limits_{\omega_{s,d}}W_{t_{ev}(p,r,m')-v_2(s)}(\ZZ_2^{cycl})   & \textrm{if }*=2r\\
\prod\limits_{m'\in \ZZ\setminus 2\ZZ}\prod\limits_{s\mid m'}\prod\limits_{\nu=0}^{t_{ev}-1}\prod\limits_{\omega_{s,d}}\ZZ_2^{cycl}   &   \textrm{if }*=2r+1
\end{cases}.
\]
\end{cor}

This argument will work for any ring $\mathcal{O}_F$, where $F/\QQ_p$ is an algebraic extension such that $\QQ_p(\zeta_{p^\infty})\subseteq F$. In fact, this will work for any ring $R$ such that $R^\wedge_p$ is perfectoid and $R$ has bounded $p^\infty$-torsion.

\subsection{Relationship to the work of Speirs}\label{sec: Speirs}
There are many connections between this paper and the papers \cite{Speirs_coordinate_axes} and \cite{Speirs_truncated_polynomials} of Speirs. Many of the arguments found here were inspired from those found in these papers. That being said, we can actually recover many of the results in both of these papers from our result, either as a direct application or as a consequence of the proof.

The first connection is to \cite{Speirs_truncated_polynomials}, which itself is a revisit of \cite{HM_cyclic_polytopes}. We may recover the special case of the dual numbers from \cite{Speirs_truncated_polynomials}. Directly applying our result when $d=1$, we see that \[\pi_*\left(K(R[x]/x^2, x)^\wedge_p\right)\cong \prod_{m'\in J_p} W_{h}(R)\] where \[h=\begin{cases}
t_{od}   &   \textrm{if }*,p\textrm{ are odd}\\
1   &   \textrm{if }*\geq 1\textrm{ is odd, }p=2\\
0   &   \textrm{otherwise}
\end{cases}.\] This is because when $d=1$, $\omega_{s,d}=0$ unless $s=1$. Thus the inner product is only non-zero when $s=1$, and hence only when $m'$ is odd. Note that this agrees with the $h$ function in \cite{Speirs_truncated_polynomials}, and so we have as a consequence\[\pi_{2r-1}\left(K(k[x]/x^2, x)^\wedge_p\right)\cong \mathbb{W}_{2r}(k)/V_2\mathbb{W}_r(k)\] by \parencite[Lemma 2]{Speirs_truncated_polynomials}, where $k$ a perfect field of characteristic $p>0$, and $\mathbb{W}$ the big Witt vectors.

What is perhaps more surprising is that we can also recover and extend calculations from \cite{Speirs_coordinate_axes} as well. To see this, let $\Pi_S$ be the pointed monoid $\{0, 1, x_1,\ldots, x_d, x_1^2, x_2^2,\ldots\}$, in other words, the pointed monoid considered in \cite{Speirs_coordinate_axes} to study $\thh$ of the coordinate axes. We then get a map $B^{cy}(\Pi_S)\to B^{cy}(\Pi)$, which breaks up along the cyclic word decomposition on both sides. For a given cyclic word $\overline{\omega}$ without cyclic repetitions, the map $B^{cy}(\Pi_S, \overline{\omega})\to B^{cy}(\Pi, \overline{\omega})$ is a $\mathbb{T}$-equivariant homeomorphism. On the other hand, if $\overline{\omega}$ has cyclic repetitions, then by \parencite[Lemma 9(2)]{Speirs_coordinate_axes} $B^{cy}(\Pi_S,\overline{\omega})$ is $\mathbb{T}$-equivariantly contractable.

In summation, the map $B^{cy}(\Pi_S)\to B^{cy}(\Pi)$ is exactly the inclusion of the summand \[\bigvee_{\substack{\overline{\omega}\\ \textrm{No cyclic repetitions}}}B^{cy}(\Pi, \overline{\omega})\hookrightarrow \bigvee_{\overline{\omega}}B^{cy}(\Pi, \overline{\omega})\] Following this summand through the argument above then gives the following.

\begin{cor}
Let $p$ be an odd prime and $R$ a perfectoid ring. Define \[C_d:= R[x_1,\ldots, x_d]/(x_ix_j)_{i\neq j}.\] Then there are isomorphisms
\[
\pi_*\left(K(C_d, \mathfrak{m})^\wedge_p\right)\cong \begin{cases}
\prod\limits_{m'\in 2J_p}\prod\limits_{\substack{s\mid m'p^{t_{ev}-1}\\  s\textrm{ even}}}\prod\limits_{cyc_d(s)}W_{t_{ev}(p,r,m,)-v_p(s)}(R)   & \textrm{if }*=2r\\
\prod\limits_{m'\in J_p\setminus 2J_p}\prod\limits_{s\mid m'p^{t_{od}-1}}\prod\limits_{cyc_d(s)}W_{t_{od}(p, r, m')-v_p(s)}(R)   &   \textrm{if }*=2r+1
\end{cases}.
\]

For p=2, there are isomorphisms
\[
\pi_*\left(K(C_d,\mathfrak{m})^\wedge_2\right)\cong \begin{cases}
\prod\limits_{m'\in 2\ZZ}\prod\limits_{\substack{s\mid m'2^{t_{ev}-1}\\  s\textrm{ even}}}\prod\limits_{cyc_d(s)}W_{t_{ev}(2,r,m')-v_2(s)}(R)   & \textrm{if }*=2r\\
\prod\limits_{m'\in \ZZ\setminus 2\ZZ}\prod\limits_{s\mid m'}\prod\limits_{\nu=0}^{t_{ev}-1}\prod\limits_{cyc_d(s)} R   &   \textrm{if }*=2r+1
\end{cases}.
\]
\end{cor}
In particular, when $R$ is a perfect field of characteristic $p$ we recover \parencite[Theorem 1]{Speirs_coordinate_axes}.

\printbibliography
\end{document}